\documentclass[journal,onecolumn]{IEEEtran}
\usepackage{color}
\usepackage{graphicx}
\usepackage[dvips]{psfrag}
\usepackage{psfrag}
\usepackage{subfigure}
\usepackage{stfloats}
\usepackage{amsfonts}
\usepackage{amsbsy}
\usepackage{amsmath}
\usepackage{amssymb}
\usepackage{latexsym}
\usepackage{url}
\usepackage{multirow}
\usepackage{rotating}
\usepackage{cite}
\usepackage{dsfont}
\usepackage{widetext}

\begin{document}
\title{Distributed Approach for DC Optimal Power Flow Calculations}
%\vspace{5cm}
\author{Javad Mohammadi, \IEEEmembership{IEEE Student Member}, Soummya Kar, \IEEEmembership{IEEE Member}, Gabriela Hug, \IEEEmembership{IEEE Member}%\vspace{-0.5cm}% <-this % stops a space
\thanks{The authors are with the Department of Electrical and Computer Engineering, Carnegie Mellon University, Pittsburgh, PA, 15213, USA, e-mails: jmohamma@andrew.cmu.edu, ghug@ece.cmu.edu, soummyak@andrew.cmu.edu}}% <-this % stops a space


\maketitle

\begin{abstract}
%The trend in the electric power system is to move towards increased amounts of distributed resources in combination with adjustable demand which suggests a transition from the current highly centralized to a more distributed control structure. In this paper, we propose a method which enables a fully distributed solution of the DC Optimal Power Flow problem, i.e.~the generation settings which minimize cost while supplying the load and ensuring that all line flows are below their limits are determined in a distributed fashion. The approach consists of a distributed procedure that aims at solving the first order optimality conditions in which individual bus optimization variables are iteratively updated through simple local computations and information is exchanged with neighboring entities. In particular, the update for a specific bus consists of a term which takes into account the coupling between the neighboring Lagrange multiplier variables and a local innovation term that enforces the demand/supply balance. The buses exchange information on the current update of their multipliers and the bus angle with their neighboring buses. There is no need to exchange the generation settings or the cost parameters of the generators thereby providing some inherent privacy with respect to these values. An analytical proof is given that the proposed method converges to the optimal solution of the DC Optimal Power Flow problem. In addition, the performance is evaluated using the IEEE Reliability Test System as a test case.
The trend in the electric power system is to move towards increased amounts of distributed resources which suggests a transition from the current highly centralized to a more distributed control structure. In this paper, we propose a method which enables a fully distributed solution of the DC Optimal Power Flow problem (DC-OPF), i.e. the generation settings which minimize cost while supplying the load and ensuring that all line flows are below their limits are determined in a distributed fashion. The approach consists of a distributed procedure that aims at solving the first order optimality conditions in which individual bus optimization variables are iteratively updated through simple local computations and information is exchanged with neighboring entities. In particular, the update for a specific bus consists of a term which takes into account the coupling between the neighboring Lagrange multiplier variables and a local innovation term that enforces the demand/supply balance. The buses exchange information on the current update of their multipliers and the bus angle with their neighboring buses. An analytical proof is given that the proposed method converges to the optimal solution of the DC-OPF. Also, the performance is evaluated using the IEEE Reliability Test System as a test case.
\end{abstract}

\vspace{0.15cm}
%% Note that keywords are not normally used for peerreview papers.

\IEEEpeerreviewmaketitle

\section{Introduction}
The control responsibility of the electric power system is shared among many control entities, each responsible for a specific part of the system. While these control areas are coordinated to a certain degree, the coordination generally does not lead to system wide optimal performance, i.e.~only suboptimal solutions are achieved. Within each control area, a highly centralized control structure is used to determine the settings of the controllable devices in that area usually taking the neighboring control areas into account as static power injections.

The recent interest in distributed methods to solve economic dispatch and optimal power flow problems stem mostly from the fact that the amount of distributed generation and intelligent and adjustable demand is expected to increase significantly and ways need to be found to handle the increasing number of control variables even within a single control area. While this is also the main motivation for this paper, the same methods can also be employed to achieve improved coordination among control areas thereby leading to optimal overall system utilization.

In this paper, we present an approach which enables a distributed solution of the DC Optimal Power Flow problem. Hence, the objective is to minimize the generation cost to fully supply the load while ensuring that no line limits are violated. The proposed approach is based on obtaining a solution to the first order optimality conditions of the corresponding optimization problem in a fully distributed fashion. These conditions include constraints which constitute a coupling of the Lagrange multipliers associated with the power flow equations and line constraints at neighboring buses and lines. This is used to formulate a term in the updates of the local variables and multipliers which takes into account these couplings. In addition, the power flow equations at the buses are used to form another term which corresponds to an innovation term enforcing the demand/supply balance. The information the buses share with neighboring buses is limited to the updates of the bus angle and the local Lagrange multipliers, i.e.~there is no need to share information about the generation settings or the cost parameters during the iterative process.

The paper is organized as follows: Sect.~\ref{relatedWork} provides an overview over related work. Sect.~\ref{DCOPF} gives the DC Optimal Power Flow formulation and the resulting first-order optimality conditions which are used in Sect.~\ref{Distributed} to derive the proposed distributed approach. Section \ref{Proof} provides the proof of convergence for the proposed algorithm. Simulation results are given in Sect.~\ref{Simulations} and Sect.~\ref{Conclusion} concludes the paper.

\section{Related Work}\label{relatedWork}
There has been a range of publications on the usage of consensus based approaches to solve the economic dispatch problem including \cite{zhang_incremental_2011,zhang_leader_2011,zhang_decentralizing_2011,kar_distributed_2012,zhang_convergence_2012,yang_consensus_2013,Davoudi_EDloss_2014}. These approaches are based on the fact that if the grid is neglected and only the solution to the economic dispatch problem is sought, then the optimal solution is obtained if the marginal cost of all the generators are equal to each other. Consequently, a consensus approach can be employed to seek an agreement for the marginal cost value. The additional constraint of total generation having to be equal to total load is taken into account differently in the various approaches. In \cite{zhang_incremental_2011,zhang_leader_2011,zhang_decentralizing_2011,zhang_convergence_2012} a leading role is assigned to one of the distributed agents whereas in \cite{kar_distributed_2012,yang_consensus_2013} local innovation gradients computed solely on the basis of local demand/supply information are used to enforce that constraint. In our paper, we incorporate power flow equations and limits on transmission lines. As this results in non-equal values for the marginal costs of the generators in the optimum for cases with congestions, a direct application of the consensus approach is not sufficient any more.

Distributed approaches to solve the Optimal Power Flow problem to determine the optimal generation settings taking into account grid constraints have mostly been based on decomposition theory such as Lagrangian Relaxation and Augmented Lagrangian Relaxation \cite{conejo_decomposition_2006}. Early examples for such applications include \cite{conejo_multi-area_1998,kim_comparison_2000,hur_evaluation_2002,nogales_decomposition_2003,marvasti2014optimal}. A more recent approach is presented in \cite{kraning_dynamic_2013} where an alternating direction method of multipliers, an augmented lagrangian relaxation method, is employed to solve a multi-step DC optimal power flow problem. The proposed approach is inherently different from these decomposition theory based approaches in multiple ways: methodologically, it is based on directly solving the first order optimality conditions and, hence, technically involves (and reduces the original optimization motive to) obtaining solutions of a coupled system of linear equations with geometric constraints in a fully distributed manner.

\section{DC Optimal Power Flow}\label{DCOPF}
The goal in DC Optimal Power Flow is to determine the generation dispatch which minimizes the cost to supply a given load taking into account operational constraints such as line limits and generation capacities. The grid is modeled using a DC approximation, hence, it is assumed that angle differences across lines are small, voltage magnitudes are all equal to 1pu and resistances of the lines are negligible.

Modeling generation costs using a quadratic cost function, the mathematical problem formulation results in
\begin{equation}
\min_{P_G} \sum_{n \in \Omega_G} \left( a_n P_{G_n}^2 + b_n P_{G_n} + c_n \right)\label{OPFobj}
\end{equation}
\begin{eqnarray}
\textrm{s.t.}~\sum_{n\in\Omega_{G_i}}P_{G_n} - P_{L_i} \hspace{-0.6cm} &=& \hspace{-0.6cm}\sum_{j \in \Omega_i} \frac{\theta_i - \theta_j}{X_{ij}},~ \forall i\in \{1,\ldots, N_B\}\label{OPFcons0}\\
\theta_1 \hspace{-0.4cm} &=& \hspace{-0.4cm} 0 \\
\underline{P}_{G_n} \leq & P_{G_n} & \leq \overline{P}_{G_n}, ~~\forall n \in \Omega_G\\
-\overline{P}_{ij} \leq & \frac{\theta_i - \theta_j}{X_{ij}} & \leq \overline{P}_{ij}, ~~~\forall ij \in \Omega_L\label{OPFcons}
\end{eqnarray}
where the variables have the following meanings:\\
\begin{tabular}{ll}
$P_{L_i}$: & load at bus $i$\\
$P_{G_n}$: & output of generator $n$\\
$a_n, b_n, c_n$: & cost parameters of generator $n$\\
$\theta_i$: & angle at bus $i$\\
$\underline{P}_{G_n}, \overline{P}_{G_n}$: & lower and upper limits on generation\\
%\end{tabular}\\
%\begin{tabular}{ll}
$\Omega_G$: & set of all generators\\
$\Omega_{G_i}:$ & set of generators at bus $i$\\
$\Omega_L$: & set of lines in the grid \\
$\Omega_i$: & set of buses connected to bus $i$\\
$X_{ij}$: & reactance of line connecting buses $i$ and $j$\\
$\overline{P}_{ij}$: & capacity of line connecting buses $i$ and $j$\\
\end{tabular}\\
and $i=1$ is taken to be the slack bus.

The Lagrange function for this optimization problem is given by
\begin{eqnarray}
L \hspace{-0.2cm} &=& \hspace{-0.2cm} \sum_{n \in \Omega_G} \left( a_n P_{G_n}^2 + b_n P_{G_n} + c_n \right) \nonumber\\
&& \hspace{-1cm} + \sum_{n \in \Omega_G} \mu_{n}^+ \left( P_{G_n} - \overline{P}_{G_n} \right) + \sum_{n \in \Omega_G} \mu_{n}^- \left( - P_{G_n} + \underline{P}_{G_n} \right)\nonumber\\
&& \hspace{-1cm}  + \sum_{i=1}^{N_B} \lambda_i \left(- \hspace{-0.2cm} \sum_{n\in\Omega_{G_i}}\hspace{-0.2cm} P_{G_n} + P_{L_i} + \sum_{j \in \Omega_i} \frac{\theta_i - \theta_j}{X_{ij}} \right) + \lambda_0 \theta_1 \nonumber\\
&& \hspace{-1cm}  + \hspace{-0.2cm} \sum_{ij \in \Omega_L}\hspace{-0.2cm} \mu_{ij} \left(\frac{\theta_i-\theta_j}{X_{ij}} -\overline{P}_{ij}\right) + \hspace{-0.2cm}\sum_{ij \in \Omega_L} \hspace{-0.2cm}\mu_{ji} \left(-\frac{\theta_i-\theta_j}{X_{ij}} -\overline{P}_{ij}\right)
\end{eqnarray}
where $\lambda$'s and $\mu$'s correspond to Lagrange multipliers. Hence, the first order optimality conditions result in
\begin{alignat}{6}
\frac{\partial L}{\partial P_{G_n}} &=& ~~2a_n P_{G_n} + b_n - \lambda_n + \mu_n^+ - \mu_n^- ~&=&~ 0 \label{KKT1}\\
\frac{\partial L}{\partial \theta_i}&=& ~~\lambda_i \cdot \sum_{j\in \Omega_i} \frac{1}{X_{ij}} - \sum_{j\in \Omega_i} \lambda_j \frac{1}{X_{ij}} &&\nonumber \\
&& + \sum_{j \in \Omega_{i}} \left( \mu_{ij} - \mu_{ji} \right) \frac{1}{X_{ij}} ~&=&~ 0 \label{KKT4}\\
\frac{\partial L}{\partial \lambda_i} &=& - \hspace{-0.2cm} \sum_{n\in\Omega_{G_i}}\hspace{-0.2cm}P_{G_n} + P_{L_i} + \sum_{j\in \Omega_i} \frac{\theta_i-\theta_j}{X_{ij}}  &=& 0\label{KKT5}\\
\frac{\partial L}{\partial \lambda_0} &=& \theta_1  &=& 0\label{KKT8}\\
\frac{\partial L}{\mu_{n}^+} &=& P_{G_n} - \overline{P}_{G_n}  ~&\leq&~ 0 \label{KKT2}\\
\frac{\partial L}{\mu_{n}^-} &=& -P_{G_n} + \underline{P}_{G_n}  ~&\leq&~ 0 \label{KKT3}\\
\frac{\partial L}{\partial \mu_{ij}} &=& \frac{\theta_i-\theta_j}{X_{ij}} -\overline{P}_{ij} &\leq& 0 \label{KKT6}\\
\frac{\partial L}{\partial \mu_{ji}} &=&  -\frac{\theta_i-\theta_j}{X_{ij}} -\overline{P}_{ij} &\leq& 0 \label{KKT7}
\end{alignat}
for all $i\in\{1,\ldots,N_B\}$, $n\in \Omega_G$ and $ij \in \Omega_L$ plus the complementary slackness conditions for the inequality constraints and the positivity constraints on the $\mu$'s. As $\lambda_0=0$ due to the fact that the choice of the slack bus does not have any influence on the result, it is omitted in (\ref{KKT4}) for the slack bus. Consequently, in order to find a solution to the DC-OPF problem, the above constrained equation system needs to be solved.

\section{Distributed Approach}\label{Distributed}
%\textbf{SOUMMYA: could you write something general about the innovations approach?}
This section presents a distributed iterative approach to solving the first order constrained equation system given in Sect.~\ref{DCOPF}, where each bus merely exchanges information with its physically connected neighbors during the course of iterations.
\\ \indent In the proposed approach, each bus $i$ updates the variables $\lambda_i,~\theta_i$ and $P_{G_n}, n \in \Omega_{G_i}$ which are directly associated with that bus and the $\mu_{ij}$'s which correspond to the constraints on the flows into bus $i$ from lines $ij$. Denoting the iteration counter by $k$ and the iterates by $x_{i}(k)$ which include the variables associated with bus $i$ at iteration $k$, i.e. $x_i(k)=[\lambda_{i}(k),\theta_{i}(k),\mu_{ij}(k),P_{G_i}(k)]$, the general format of the local updates is given by
\begin{equation}
x_{i}(k+1) = \mathbb{P}_i\left (x_{i}(k)+\Phi_i g_i(x_{j}(k))\right )~~~j\in\Omega_{i}
\label{update}
\end{equation}
In the above, the function $g_{i}(\cdot)$ represents the first order optimality constraints related to bus $i$. Also, $\Phi_i$ is the vector of tuning parameters. Moreover, $\mathbb{P}_i$ is the projection operator which projects $x_i$ onto its determined feasible space.
\\ \indent Note that, $g_{i}(x_{j}(k))$ depends only on the iterates $x_{j}(k)$ of neighboring buses $j$ in the physical neighborhood of $i$.
Hence, a distributed implementation of (\ref{update}) is possible. As will be seen, $\mu_n^+,~\mu_n^-$ do not need to be known and therefore no update is needed. In addition, $\lambda_0$ will always be equal to zero, i.e.~does not need to be considered in the updates, neither.

The Lagrange multipliers $\lambda_i$ are updated according to
\begin{eqnarray}
%\lambda_i(k+1) \hspace{-0.2cm} &=& \hspace{-0.2cm} \lambda_i(k) - \beta \cdot \left( \lambda_i(k) \sum_{j\in \Omega_i} \frac{1}{X_{ij}} - \sum_{j\in \Omega_i} \lambda_j(k) \frac{1}{X_{ij}} \right. \nonumber \\
%&& \hspace{-1.2cm} \left. + \sum_{j \in \Omega_{i}} \left( \mu_{ij}(k) - \mu_{ji}(k) \right) \frac{1}{X_{ij}} \right)\nonumber \\
%&& \hspace{-1.2cm} -\alpha \cdot \left(\hspace{-0.1cm} \sum_{n\in\Omega_{G_i}}\hspace{-0.2cm}P_{G_n}(k) - P_{L_i} - \sum_{j\in \Omega_i} \frac{\theta_i(k)-\theta_j(k)}{X_{ij}} \right) \label{lambdaUpdate}
\lambda_i(k+1) \hspace{-0.2cm} &=& \hspace{-0.2cm} \lambda_i(k) - \beta \cdot \left( \frac{\partial L}{\partial \theta_i}\right)+\alpha\cdot\left(\frac{\partial L}{\partial \lambda_i} \right )\nonumber \\
&=& \hspace{-0.2cm} \lambda_i(k) - \beta \cdot \left( \lambda_i(k) \sum_{j\in \Omega_i} \frac{1}{X_{ij}} - \sum_{j\in \Omega_i} \lambda_j(k) \frac{1}{X_{ij}} \right. \nonumber \\
&& \hspace{-1.2cm} \left. + \sum_{j \in \Omega_{i}} \left( \mu_{ij}(k) - \mu_{ji}(k) \right) \frac{1}{X_{ij}} \right)\nonumber \\
&& \hspace{-1.2cm} -\alpha \cdot \left(\hspace{-0.1cm} \sum_{n\in\Omega_{G_i}}\hspace{-0.2cm}P_{G_n}(k) - P_{L_i} - \sum_{j\in \Omega_i} \frac{\theta_i(k)-\theta_j(k)}{X_{ij}} \right) \label{lambdaUpdate}
\end{eqnarray}
where $\alpha, \beta >0 $ are tuning parameters and $k$ denotes the iteration index. Hence, the first term corresponds to the optimality condition (\ref{KKT4}) which reflects the coupling between the Lagrange multipliers and the second term constitutes an innovation term based on the power balance equations (\ref{KKT5}). The update makes intuitive sense, e.g.~if the power balance (\ref{KKT5}) is not fulfilled because generation is too high, it leads to a reduction in $\lambda_i$ which on the other hand, as shown next, leads to a decrease in the $P_{G_n}, n \in \Omega_{G_i}$. Furthermore, if no line constraints are binding the $\mu$'s are equal to zero and the part of the update in the first row leads to finding an agreement between the $\lambda$'s at all buses.

Knowing the value of the Lagrange multiplier $\lambda_i$, the following update for the generators $P_{G_n}, n\in\Omega_{G_i}$ can be carried out:
\begin{equation}
P_{G_n}(k+1) = \mathbb{P}_n\hspace{-0.1cm}\left(\hspace{-0.05cm}P_{G_n}(k)-\hspace{-0.1cm}\frac{1}{2a_n}\cdot\frac{\partial L}{\partial P_{G_n}}\right )\hspace{-0.05cm}=\hspace{-0.05cm}\mathbb{P}_n\hspace{-0.05cm}\left(\hspace{-0.05cm}\frac{\lambda_i(k) - b_n}{2a_n}\right ) \label{PGUpdate}
\end{equation}
Here $\mathbb{P}_n$ is the operator which projects the value determined by (\ref{PGUpdate}) into the feasible space defined by the upper and lower limits $\overline{P}_{G_n}$ and $\underline{P}_{G_n}$, i.e.~if the value is greater than $\overline{P}_{G_n}$, $P_{G_n}(k+1)$ is set to that upper limit and similarly for the lower limits.
%The upper and lower limits $\overline{P}_{G_n}$ and $\underline{P}_{G_n}$ are taken into account by projecting the value determined by (\ref{PGUpdate}) into the feasible space, i.e.~if the value is greater than $\overline{P}_{G_n}$, $P_{G_n}(k+1)$ is set to that upper limit and similarly for the lower limits.
This is equivalent to using the full equation (\ref{KKT1}) including the multipliers $\mu_n^+$ and $\mu_n^-$ to update $P_{G_n}$. As these multipliers do not appear in any other constraint it is not necessary to provide an update for them.

The bus angles are updated according to
\begin{eqnarray}
\theta_i(k+1)\hspace{-0.3cm} &=&\hspace{-0.27cm}\theta_i(k)-\gamma \left( \frac{\partial L}{\partial \lambda_i}\right)\nonumber \\
 &=&\hspace{-0.27cm}\theta_i(k) - \hspace{-0.03cm}\gamma\hspace{-0.08cm} \left(\hspace{-0.09cm} -\hspace{-0.3cm} \sum_{n\in\Omega_{G_i}}\hspace{-0.2cm}P_{G_n}(k) + \hspace{-0.07cm}P_{L_i} \hspace{-0.07cm}+ \hspace{-0.16cm}\sum_{j\in \Omega_i} \hspace{-0.06cm} \frac{\theta_i(k)-\theta_j(k)}{X_{ij}} \hspace{-0.1cm}\right)\label{thetaUpdate}\nonumber\\
%\theta_i(k+1) = \theta_i(k) -\gamma \hspace{-0.3cm}&(  -\sum_{n\in\Omega_{G_i}}\hspace{-0.2cm}P_{G_n}(k) + P_{L_i} \\ \nonumber
%&+ \sum_{j\in \Omega_i} \frac{\theta_i(k)-\theta_j(k)}{X_{ij}} )&\label{thetaUpdate}
\end{eqnarray}
with $\gamma >0$ being a tuning parameter. Hence, the power balance equation (\ref{KKT5}) is used for the update. It again makes intuitive sense because if the power balance is not fulfilled and the load plus what is flowing onto the lines is greater than the generation at that bus, the angle is reduced which results in a reduction of the residual of that constraint.

The Lagrange multipliers $\mu_{ij},~\mu_{ji}$ appear in the $\lambda$ updates (\ref{lambdaUpdate}), and hence, values and updates for these multipliers are needed. The update is given by
\begin{eqnarray}
\mu_{ij}(k+1)\hspace{-0.2cm} &=& \hspace{-0.25cm}\mathbb{P}\hspace{-0.1cm}\left(\hspace{-0.05cm} \mu_{ij}(k) + \delta \hspace{-0.05cm} \left(\frac{\partial L}{\partial \mu_{ij}} \right)\hspace{-0.1cm}\right)\nonumber\\
\mu_{ij}(k+1)\hspace{-0.2cm} &=& \hspace{-0.25cm}\mathbb{P}\hspace{-0.1cm}\left(\hspace{-0.05cm}\mu_{ij}(k) - \delta  \hspace{-0.05cm} \left(\overline{P}_{ij} - \frac{\theta_i(k) - \theta_j(k)}{X_{ij}} \right)\hspace{-0.1cm}\right)\label{muUpdate1}\nonumber\\
\\
\mu_{ji}(k+1)\hspace{-0.2cm} &=& \hspace{-0.25cm}\mathbb{P}\hspace{-0.1cm}\left(\hspace{-0.05cm} \mu_{ji}(k) + \delta \hspace{-0.05cm} \left(\frac{\partial L}{\partial \mu_{ji}} \right)\hspace{-0.1cm}\right)\nonumber\\
\mu_{ji}(k+1)\hspace{-0.2cm} &=& \hspace{-0.25cm}\mathbb{P}\hspace{-0.1cm}\left(\hspace{-0.05cm} \mu_{ji}(k) - \delta  \hspace{-0.05cm} \left(\overline{P}_{ji} + \frac{\theta_i(k) - \theta_j(k)}{X_{ij}} \right)\hspace{-0.1cm}\right)\label{muUpdate2}\nonumber\\
\end{eqnarray}
with $\delta > 0$ being a tuning parameter. Consequently, the inequalities (\ref{KKT6}) and (\ref{KKT7}) are used. The projection operator ($\mathbb{P}$) enforces the positivity constraint on the $\mu$'s by setting the $\mu_{ij}(k+1)$ and $\mu_{ji}(k+1)$ equal to zero if the update (\ref{muUpdate1}) and (\ref{muUpdate2}) yield negative values, respectively. Assuming that the current value for the line flow $P_{ij} = (\theta_i - \theta_j)/X_{ij}$ from bus $i$ to bus $j$ is positive but below its limit $\overline{P}_{ij}$ the update (\ref{muUpdate1}) yields a decreasing value for $\mu_{ij}$ with a minimum value of zero due to the projection into the feasible space $\mu_{ij}\geq 0$. If the flow is above the line limit, the value for $\mu_{ij}$ increases indicating a binding constraint.

It should be noted that all of these updates have purposely been defined only based on the variables from the previous iteration in order to allow for a parallel computation of all of the updates. If implemented in series, i.e.~(\ref{thetaUpdate}) uses the already updated generation values, the number of iteration until convergence decreases but computation time increases because all the computations at a specific bus have to be done after each other.

%Consequently, the update rules for the variables at bus $i$ can be written in a dense form as
%\begin{equation}
%x_i(k+1) = \tilde{x}_i(k) - A_i \cdot \tilde{x}_{ij}(k) + C_i
%\end{equation}
%where $x_i(k+1)$ is the vector of updated variables at bus $i$ without projection, i.e.
%\begin{equation}
%x_i(k+1) = [\lambda_i(k+1),\theta_i(k+1),\mu_{ij}(k+1),P_{G_i}(k+1)]^T\nonumber
%\end{equation}
%and $\tilde{x}_i(k)$ is the projection of $x_i(k)$ into the feasible space for the generation settings and the Lagrange Multipliers for the line constraints. The vector $x_{ij}(k)$ includes the vector $x_i(k)$ and the variable vectors of the neighboring buses $j \in \Omega_i$.
%
%For bus $i$ with a single neighboring bus, the matrix $A_i$ is therefore given by (\ref{updateMatrix}) and $C_i$ by
%\begin{equation}
%C_i = \begin{tabular}{cccc}[$\alpha P_{L_i}$ & $-\gamma P_{L_i}$ & $-\delta \overline{P}_{ij}$ & $-\frac{b_i}{2a_i} + \alpha P_{L_i}$ ]$^T$\end{tabular}
%\end{equation}

Consequently, the update rules for the all variables can be written in a dense form as
\begin{align}
X(k+1)&=\widetilde{X}(k)-A\widetilde{X}(k)+C\nonumber\\
\widetilde{X}(k+1)&=\mathbb{P}(X(k+1))\label{iterative0}
\end{align}
where $X$ is the vector of the stacked variables ($\lambda_i,\theta_i,\mu_{ij},P_{G_i}$) for all buses $i = \{1,\ldots, N_B$), $j\in\Omega_i$ and $\mathbb{P}$ is the projection operator which ensures that the Lagrange Multipliers for the line constraints stay positive and the generation outputs stay within the given bound. Hence, $\widetilde{X}$ is the vector of the stacked projected variables. Equation (\ref{updateMatrix}) presents (\ref{iterative0}) in more detail.

\newcounter{tempequationcounter}
\begin{figure*}[!t]
\normalsize
\setcounter{tempequationcounter}{\value{equation}}
\begin{equation}\label{updateMatrix}
X(k+1) =  \left(I-\left[ \begin{tabular}{cccccccc}
$-\alpha $ & $\beta $ & 0 & 0 \\
$\gamma $ & 0 & 0 & 0 \\
0 & 0 & $-\delta $ & 0 \\
0 & 0 & 0 & $\frac{1}{2a}$
\end{tabular} \right] \right)\left[ \begin{tabular}{cccccccc}
\vspace{0.1cm}
$\frac{\partial L}{\partial \lambda}$ \\
\vspace{0.1cm}
$\frac{\partial L}{\partial \theta}$ \\
\vspace{0.1cm}
$\frac{\partial L}{\partial \mu}$\\
\vspace{0.1cm}
$\frac{\partial L}{\partial P_G}$
\end{tabular} \right]=\left(I-\overbrace{\left[ \begin{tabular}{cccccccc}
$\beta B $ & $-\alpha B$ & $\beta B_y^T$ & $\alpha I$ \\
0 & $\gamma B $ & 0 & $-\gamma I$ \\
0 & $-\delta B_{y}$ & 0 & 0 \\
$-\frac{I}{2a} $ & 0 & 0 & $I$
\end{tabular} \right]}^\text{A} \right)\widetilde{X}(k)+\overbrace{\left[ \begin{tabular}{cccccccc}
$\alpha P_L$ \\
$-\gamma P_L$ \\
$-\delta \overline{P}_{ij}$ \\
$- \frac{b}{2a}$
\end{tabular} \right]}^\text{C}
\end{equation}
\stepcounter{tempequationcounter}
\setcounter{equation}{\value{tempequationcounter}}
\hrulefill
\vspace*{4pt}
\end{figure*}
\noindent In (\ref{updateMatrix}), $I$ and $B$ are the identity and bus admittance matrices, respectively. Moreover, $B_{y}= H\cdot(\mathcal{I}\cdot diag\frac{1}{X_{ij}})^{T}$, where $\mathcal{I}$  is the incidence matrix, and $H=\begin{bmatrix}
I\\
-I
\end{bmatrix}$,

%denoting incidence matrix with $\mathcal{I}$, $B^{T}_{y}=\begin{bmatrix}
%\mathcal{I}\\
%-\mathcal{I}
%\end{bmatrix}\cdot diag\left \{\frac{1}{X_{ij}}  \right \}$.

%\begin{figure*}[!t]
%\normalsize
%\setcounter{tempequationcounter}{\value{equation}}
%\begin{equation}\label{updateMatrix}
%X(k+1) =  \left(I-\left[ \begin{tabular}{cccccccc}
%$\beta \sum_{j \in \Omega_{ij}} \frac{1}{X_{ij}} $ & $-\alpha \sum_{j \in \Omega_{ij}} \frac{1}{X_{ij}}$ & $\frac{\beta}{X_{ij}}$ & $\alpha$ & $-\frac{\beta}{X_{ij}}$ & $\frac{\alpha}{X_{ij}}$ & $-\frac{\beta}{X_{ij}}$ & 0\\
%0 & $\gamma \sum_{j \in \Omega_{ij}} \frac{1}{X_{ij}}$ & 0 & $-\gamma$ & 0 & $-\gamma \frac{1}{X_{ij}}$ & 0 & 0 \\
%0 & $-\delta \frac{1}{X_{ij}}$ & 0 & 0 & 0 & $\delta \frac{1}{X_{ij}}$ & 0 & 0\\
%$- \frac{1}{2a_i} + \frac{\beta}{2a_i} \sum_{j \in \Omega_{ij}} \frac{1}{X_{ij}} $ & $-\frac{\alpha}{2a_i} \sum_{j \in \Omega_ij} \frac{1}{X_{ij}}$ & $\frac{\beta}{2a_i} \frac{1}{X_{ij}}$ & $\frac{\alpha}{2a_i}+1$ & $-\frac{\beta}{2a_i}\frac{1}{X_{ij}}$ & $\frac{\alpha}{2a_i}\frac{1}{X_{ij}}$ & $-\frac{\beta}{2a_i}\frac{1}{X_{ij}}$ & 0
%\end{tabular} \right] \right)
%\end{equation}
%\stepcounter{tempequationcounter}
%\setcounter{equation}{\value{tempequationcounter}}
%\hrulefill
%\vspace*{4pt}
%\end{figure*}

\section{Convergence Analysis}\label{Proof}
This section presents a formal proof that any limit point of the proposed algorithm is the optimal solution of the OPF problem. Specifically, we first show that a fixed point of the proposed iterative scheme necessarily satisfies the optimality conditions (\ref{KKT1})--(\ref{KKT7}) of the OPF problem. This is achieved in Theorem~1 in the following.
\\ \indent \textit{Theorem~1:}  Let $X^{\ast}$ be a fixed point of the proposed algorithm defined by (\ref{iterative0}). Then, $X^{\star}$ satisfies all of the optimality conditions of the OPF problem (\ref{KKT1})--(\ref{KKT7}).
\\ \indent \textit{Proof:} To prove this theorem, we verify the claim that $X^{\star}$ fulfills all of the first order optimality conditions. Note that $X^{\star}$ is the vector of stacked variables ($\lambda_i^{\star},\theta_i^{\star},\mu_{ij}^{\star},P_i^{\star}$) for all buses $i = \{1,\ldots, N_B\}$.
\\ \indent \textit{Claim~2.1:}~$X^{\star}$ fulfills the optimality conditions which enforce the positivity of the Lagrangian multipliers associated with the line limits, i.e. $\mu_{ij}^{\star} \geq 0$.
\\ \indent \textit{Verification by contradiction:} Let us assume on the contrary that in $X^{\star}$ one of the line limit multiplier variables, say $\mu_{ij}^{\star}$, is negative. Now, note that, evaluating (\ref{muUpdate1}) at $X^{\star}$ results in a non-negative value for $\mu_{ij}$ due to the projection of $\mu_{ij}$ into the set of positive reals. In other words, we have
 \begin{equation}
 \mu_{ij}^{\star} \neq  \mathbb{P} \left( \mu_{ij}^{\star} - \delta \cdot  \left(\overline{P}_{ij} - \frac{\theta_i^{\star}(k) - \theta_j^{\star}(k)}{X_{ij}}\right)\right),\nonumber
 \end{equation}
which contradicts the fact that $X^{\ast}$ is a fixed point of (\ref{muUpdate1}).
\\ \indent \textit{Claim~2.2:}~$X^{\star}$ satisfies the optimality conditions associated with the line limit constraints, (\ref{KKT6})--(\ref{KKT7}).
\\ \indent \textit{Verification by contradiction:} Let us assume that $X^{\star}$ does not fulfill (\ref{KKT6}) for all $i$ and $j$, i.e., there exists $(i,j)$ such that $\frac{\theta_i^{\star}-\theta_j^{\star}}{X_{ij}} >\overline{P}_{ij}$. This implies that the value of the innovation term in (\ref{muUpdate1}) is negative when evaluated at $X^{\star}$.
Also, note that, based on the claim~2.1, $\mu^{\star}_{ij} \geq 0$. Therefore, evaluating (\ref{muUpdate1}) at $X^{\star}$ results in a value greater than $\mu_{ij}^{\star}$, i.e.,\\
 \begin{equation}
 \mu_{ij}^{\star} < \mathbb{P} \left( \mu_{ij}^{\star} - \delta \cdot  \left(\overline{P}_{ij} - \frac{\theta_i^{\star}(k) - \theta_j^{\star}(k)}{X_{ij}}\right)\right),\nonumber
 \end{equation}
which contradicts the fact that $X^{\ast}$ is a fixed point of (\ref{muUpdate1}). Similar arguments can be used to prove that $X^{\star}$ fulfills (\ref{KKT7}).
\\ \indent \textit{Claim~2.3:} $X^{\star}$ satisfies the optimality conditions associated with the complementary slackness condition, i.e., for all pairs $(i,j)$,\\
\begin{equation}
 \mu_{ij}^{\star}\cdot\left(\frac{\theta_i^{\star}-\theta_j^{\star}}{X_{ij}} -\overline{P}_{ij}\right)=0.\nonumber
 \end{equation}
\\ \indent \textit{Verification by contradiction:} Let us assume on the contrary that $X^{\ast}$ does not satisfy the above complementary slackness condition, i.e., there exists a pair $(i,j)$ such that both $\mu_{ij}^{\star}$ and $\frac{\theta_i^{\star}-\theta_j^{\star}}{X_{ij}} -\overline{P}_{ij}$ are non-zero. Hence, according to the claims~2.1 and~2.2, we must have, $\mu_{ij}^{\star} > 0$ and $\frac{\theta_i^{\star}-\theta_j^{\star}}{X_{ij}} <  \overline{P}_{ij}$, respectively. Now, note that evaluating (\ref{muUpdate1}) at $X^{\star}$, results in a value less than $\mu_{ij}^{\star}$, which clearly contradicts the fact that $X^{\ast}$ is a fixed point of (\ref{muUpdate1}).
\\ \indent \textit{Claim~2.4:}~$X^{\star}$ satisfies the local load balance equation (\ref{KKT5}).
\\ \indent \textit{Verification by contradiction:} Let us assume on the contrary that $X^{\star}$ does not fulfill (\ref{KKT5}), i.e., there exists $i$ such that the value of the innovation term in (\ref{thetaUpdate}) is non-zero when evaluated at $X^{\star}$. Clearly, this would lead to
\begin{equation}
\theta_i^{\star} \neq \theta_i^{\star} - \gamma \left(\hspace{-0.1cm} -\hspace{-0.2cm} \sum_{n\in\Omega_{G_i}}\hspace{-0.2cm}P_{G_n}^{\star} + P_{L_i} + \sum_{j\in \Omega_i} \frac{\theta_i^{\star}-\theta_j^{\star}}{X_{ij}} \hspace{-0.1cm}\right),\nonumber
\end{equation}
thus contradicting the fact that $X^{\ast}$ is a fixed point of (\ref{thetaUpdate}).
\\ \indent \textit{Claim~2.5:}~The coupling between the Lagrangian multipliers, given by (\ref{KKT4}), is maintained at $X^{\star}$.
\\ \indent \textit{Verification by contradiction:} Let us assume on the contrary that $X^{\star}$ does not fulfill (\ref{KKT4}) for some $i$.
%, i.e. the value of one of the innovation terms in (\ref{lambdaUpdate}) is none-zero at $X^{\star}$.
Note that (\ref{lambdaUpdate}) includes two innovation terms: the innovation term associated with the Lagrangian multipliers' coupling and the innovation term which represents the local power balance equation. We already verified that the local power balance equation is zero at $X^{\star}$ (see claim~2.4). Thus, the contradiction hypothesis necessarily implies that the innovation term associated with the Lagrangian multipliers' coupling attains a non-zero value at $X^{\star}$. This, in turn, implies that the the value of (\ref{lambdaUpdate}) is not equal to $\lambda_i^{\star}$ when evaluated at $X^{\star}$, which clearly contradicts the fact that $X^{\ast}$ is a fixed point of (\ref{lambdaUpdate}).
\\ \indent \textit{Claim~2.6:}~$X^{\star}$ satisfies the optimality conditions associated with the generation limits, (\ref{KKT2})--(\ref{KKT3}).
\\ \indent \textit{Verification by contradiction:} %Based on Theorem~2.1-2.3, $X^{\star}$ fulfills (\ref{KKT1}-\ref{KKT5},\ref{KKT6},\ref{KKT7}).
Let us assume on the contrary that there exists $i$ such that $P_{G_i}^{\star}$ does not lie in $[\underline{P}_{G_i},\overline{P}_{G_i}]$. Now, note that, plugging in $\lambda^{\star}$ in (\ref{PGUpdate}), would then result in a value different from $P_{G_i}^{\star}$, since the projection operator enforces the value of $P_{G_i}$ to stay in the specified region, $[\underline{P}_{G_i},\overline{P}_{G_i}]$. This, in turn, clearly contradicts the fact that $X^{\ast}$ is a fixed point of (\ref{PGUpdate}).

We discuss the consequences of Theorem~1. To this end, note that, since the proposed iterative scheme (\ref{iterative0}) involves continuous transformations of the iterates, it follows that, if (\ref{iterative0}) converges, the limit point is necessarily a fixed point of the iterative mapping. Since, by Theorem~1, any fixed point of (\ref{iterative0}) solves the first order optimality conditions (\ref{KKT1})--(\ref{KKT7}), we may conclude that, if (\ref{iterative0}) converges, it necessarily converges to a solution of the first order  optimality conditions (\ref{KKT1})--(\ref{KKT7}). This immediately leads to the following optimality of limit points of the proposed scheme.
 \\ \indent \textit{Theorem~2:} Suppose the OPF problem (\ref{OPFobj})--(\ref{OPFcons}) has a feasible solution that lies in the interior of the associated constraint set, and, further, assume that the proposed algorithm defined by (\ref{iterative0}) converges to a point $X^{\star}$. Then $X^{\ast}$ constitutes an optimal solution of the OPF problem (\ref{OPFobj})--(\ref{OPFcons}).
  \\ \indent \textit{Proof:} By Theorem 1 and the above remarks, $X^{\star}$ fulfills the optimality conditions (\ref{KKT1})--(\ref{KKT7}). Since the DC-OPF is a convex problem and, by assumption, is strictly feasible, it follows readily that the primal variables  $(P^{\star}$,$\theta^{\star})$ in $X^{\ast}$ constitutes an optimal solution to the OPF problem (\ref{OPFobj})--(\ref{OPFcons}).
\\ \indent In summary, we note that Theorems~1 and~2 guarantee that any fixed point of the proposed algorithm constitutes an optimal solution to the OPF problem, and, hence, in particular, if the scheme achieves convergence, the limit point is necessarily an optimal solution of the OPF problem.
\\ \indent Finally, we note, that whether the scheme converges or not depends on several design factors, in particular, the tuning parameters $\alpha$, $\beta$, $\gamma$ and $\delta$. To this end, a general sufficient condition for convergence is presented in the Appendix. Moreover, some simulation examples are presented in the following section in which we provide choices of the tuning parameters that achieve convergence (and hence, by Theorem~2, to the optimal solution of the OPF) for a class of realistic power systems.
\section{Simulation Results}\label{Simulations}
In this section, we provide simulations results to give a proof of concept of the proposed method. Note that, the tuning parameters are designed such that the
algorithm converges. The optimality of the achieved distributed solution is guaranteed according to Theorems~1 and~2.

\begin{table}[b]
\center
\begin{small}
\caption{Generator Data}\label{GenData}
\begin{tabular}{|c|c|c|c|c|c|}
\hline
 Type & Cap. & Bus & \# & a[\$/MW/MWh] & b[\$/MWh]\\
 \hline
 \#1 & 12MW & 15 & 5 & 0.36 & 20.70\\
 \hline
\#2 & 20MW & 1 & 2 & 0.45 & 21.31\\
 & & 2 & 2 & & \\
 \hline
\#3 & 50MW & 22 & 6 & 0.001 & 4\\
 \hline
\#4 & 76MW & 1 & 2 & 0.037  & 9.82\\
 &  & 2 & 2 & &\\
 \hline
\#5 &100MW & 7 & 3 & 0.027 & 17.26 \\
 \hline
\#6 &155MW & 15 & 1 & 0.0066 & 9.12\\
 &  & 16 & 1 & & \\
 &  & 23 & 2 & & \\
 \hline
\#7 & 197MW & 13 & 3 & 0.0115& 17.62\\
 \hline
\#8 & 350MW & 23 & 1 & 0.0039 & 8.95\\
 \hline
\#9 & 400MW & 18 & 1  & 0.002 & 5.17\\
 &  &  21 & 1 & &\\
\hline
\end{tabular}
\end{small}
\end{table}

\subsection{Simulation Setup}
\begin{figure}[t!]
\centering
 \footnotesize
    \includegraphics[width=7cm]{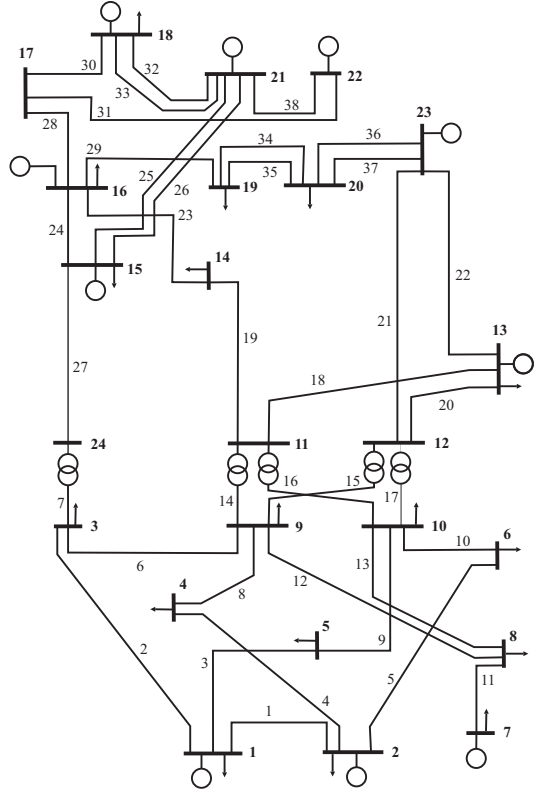}
\caption{IEEE Reliability Test System \cite{grigg_ieee_1999}.}\label{TestSystem}
\end{figure}

We use the IEEE RTS test system for our simulations \cite{grigg_ieee_1999} (Fig.~\ref{TestSystem}). The synchronous condenser is removed yielding a system with 32 generators and 17 loads. The communication network has the same topology as the physical system. The consumption of the loads is set to be equal to the values given in the original system. First, the original line limits are used which results in a situation in which no lines are at their limits for the optimal dispatch. Then, we reduce the line limits to $55\%$ of the original values for which only a suboptimal dispatch can be achieved because lines 7 and 28 reach their limits. The chosen cost parameters for the generators, their capacities, locations and how many are located at the indicated buses are given in Table \ref{GenData}. As can be seen, a mix of cheap and costly generation as well as small and large scale generation is present in the system. For both simulations, the tuning parameters are set to the values given in Table \ref{TuningValues}.

\begin{table}[t]
\center
\caption{Tuning Parameter Values}\label{TuningValues}
\begin{tabular}{|c|c|}
\hline
Parameter & Value \\
\hline
$\alpha$ & 0.1485 \\
$\beta$ & 0.0056\\
$\gamma$ & 0.005\\
$\delta$ & 0.008\\
\hline
\end{tabular}
\end{table}

We use a cold start for the simulations, i.e.~all generation values, bus angles and the Lagrange multipliers $\mu_{ij}, \mu_{ji}$ are set to zero at the start of the simulation. Merely, the Lagrange multipliers $\lambda_i$ are set to an initial value of 10~$\$/MWh$. It is reasonable to start with a non-zero value because the $\lambda$'s represent the locational marginal prices which rarely will be zero. In fact, for an actual implementation, reasonable initial settings for all of these variables could be the optimal values computed for the previous time step.
\subsection{Convergence measurements}
In order to evaluate the performance of the proposed distributed approach, two measures are introduced. The first measure determines the relative distance of the objective function from the optimal value over the iterations,
\begin{equation}
rel=\frac{\left | f-f^{*} \right |}{f^{*}},\label{rel}
\end{equation}
where $f^*$ is the optimal objective function value calculated by solving the centralized DC-OPF problem. Moreover, the value of load balance, as one of the optimality conditions, is potentially another indication of the distance from the optimal value, since the value of the load balance at the optimal point is equal to zero. Thus, we propose using the sum over the residuals of all power flow equations over the course of the iterations as the second measure of convergence, and is given by:
\begin{equation}
res=\sum_{i}\sqrt{g_{i}^{2}},
\end{equation}
\noindent where $g_{i}$ is the local power flow equation at bus $i$, which enforces supply/demand balance.

\subsection{Simulation Results}\label{SimResults}
\subsubsection{Uncongested System}
\begin{figure}[t!]
\centering
 \footnotesize
 \psfrag{iteration}[c][c]{Iterations}
 \psfrag{theta}[c][c]{$\theta$ [rad]}
 \psfrag{generation}[c][c]{$P_{G}$ [MW]}
 \psfrag{lambda}[c][c]{$\lambda$}
 \psfrag{mulines}[c][c]{$\mu_{ij},~\mu_{ji}$}
 \psfrag{0}[c][c]{}
 \psfrag{100}[c][c]{}
 \psfrag{200}[c][c]{200}
 \psfrag{300}[c][c]{}
 \psfrag{400}[c][c]{400}
 \psfrag{500}[c][c]{}
 \psfrag{600}[c][c]{600}
 \psfrag{700}[c][c]{}
 \psfrag{800}[c][c]{800}
 \psfrag{900}[c][c]{}
 \psfrag{1000}[c][c]{1000}
 \psfrag{8}[r][r]{8}
 \psfrag{12}[r][r]{12}
 \psfrag{16}[r][r]{16}
 \psfrag{20}[r][r]{20}
 \psfrag{0.2}[c][c]{0.2}
 \psfrag{0.1}[c][c]{}
 \psfrag{-0.1}[c][c]{}
 \psfrag{-0.2}[c][c]{-0.2}
 \psfrag{-0.3}[c][c]{}
 \psfrag{-0.4}[c][c]{-0.4}
 \psfrag{-0.5}[c][c]{}
 \psfrag{-0.6}[c][c]{-0.6}
 \psfrag{-0.7}[c][c]{}
 \psfrag{1}[r][r]{}
 \psfrag{2}[r][r]{200}
 \psfrag{3}[r][r]{}
 \psfrag{4}[r][r]{400}
 \psfrag{5}[r][r]{}
 \psfrag{6}[r][r]{600}
    \subfigure[]{\includegraphics[width=8cm]{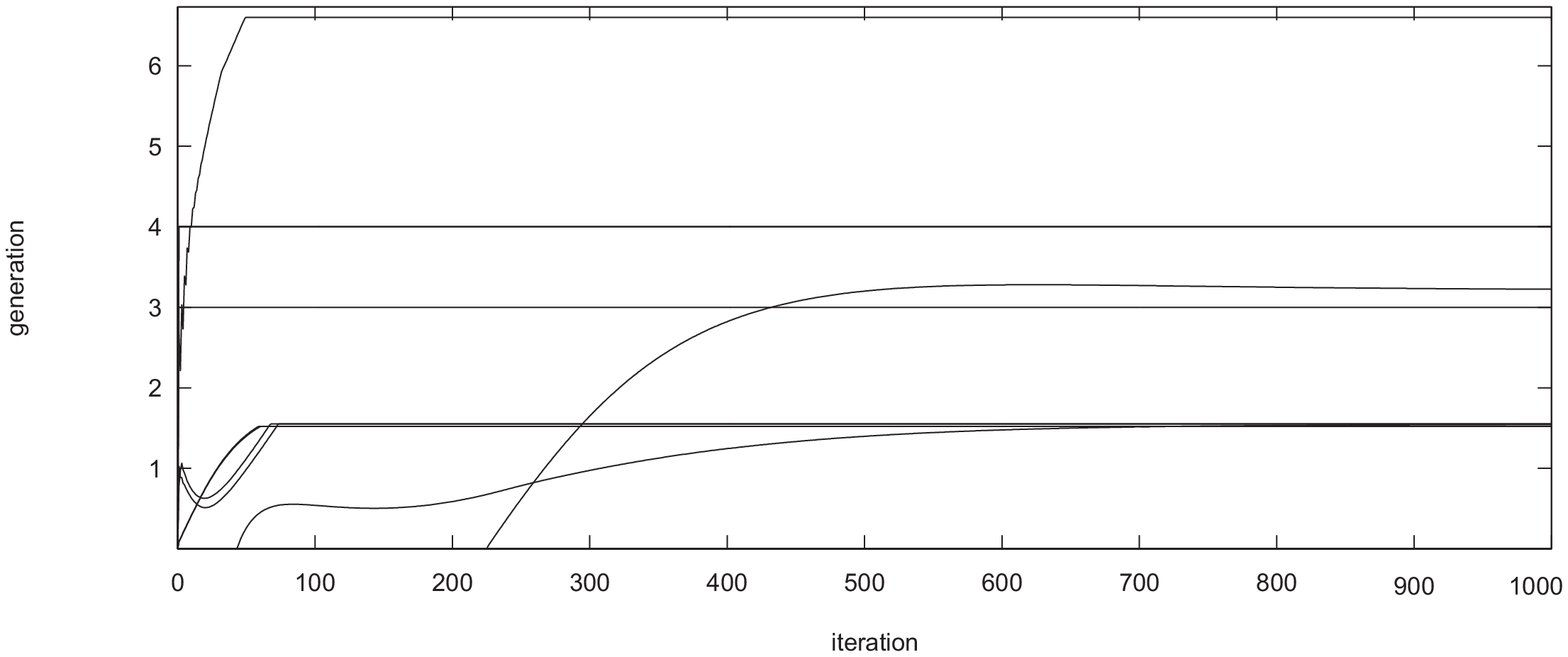}\label{generationSynch24}}
    \subfigure[]{\includegraphics[width=8cm]{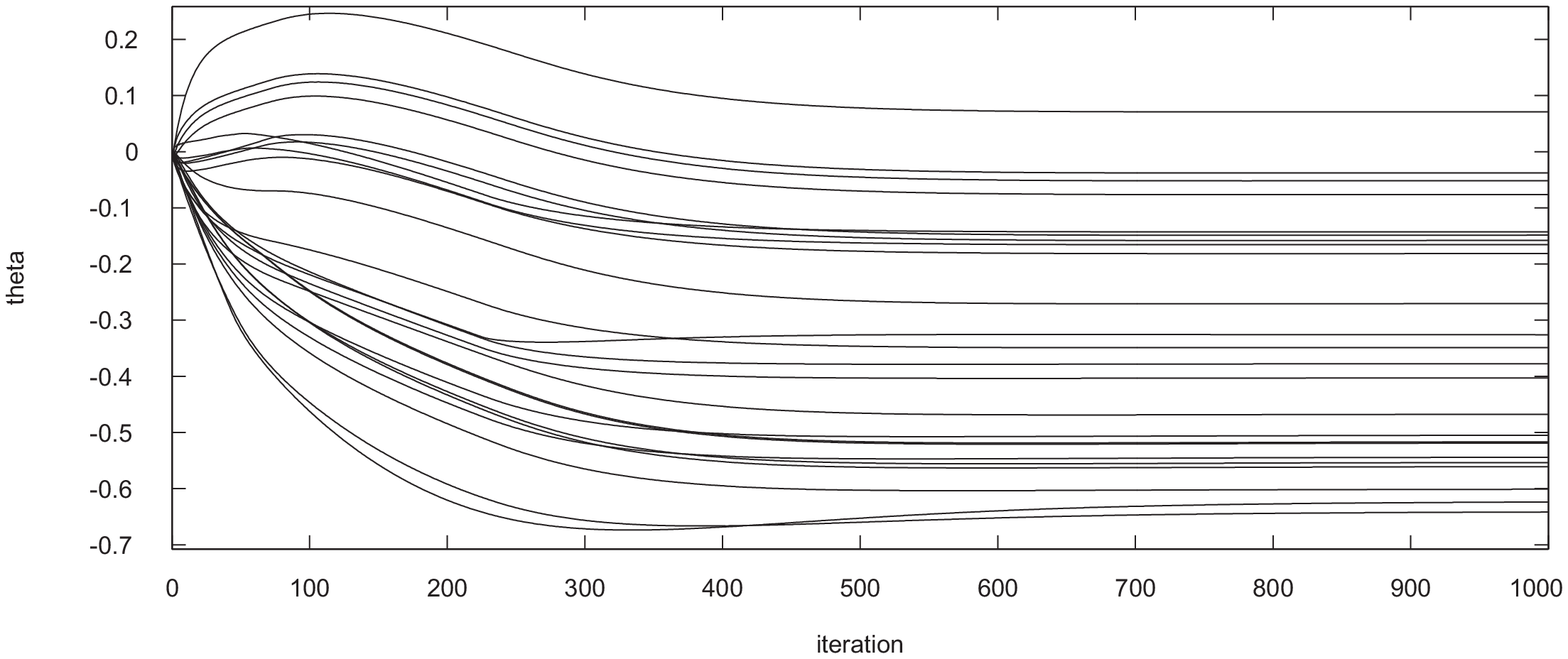}\label{thetaSynch24}}
    \psfrag{4}[r][r]{4}
    \subfigure[]{\includegraphics[width=8cm]{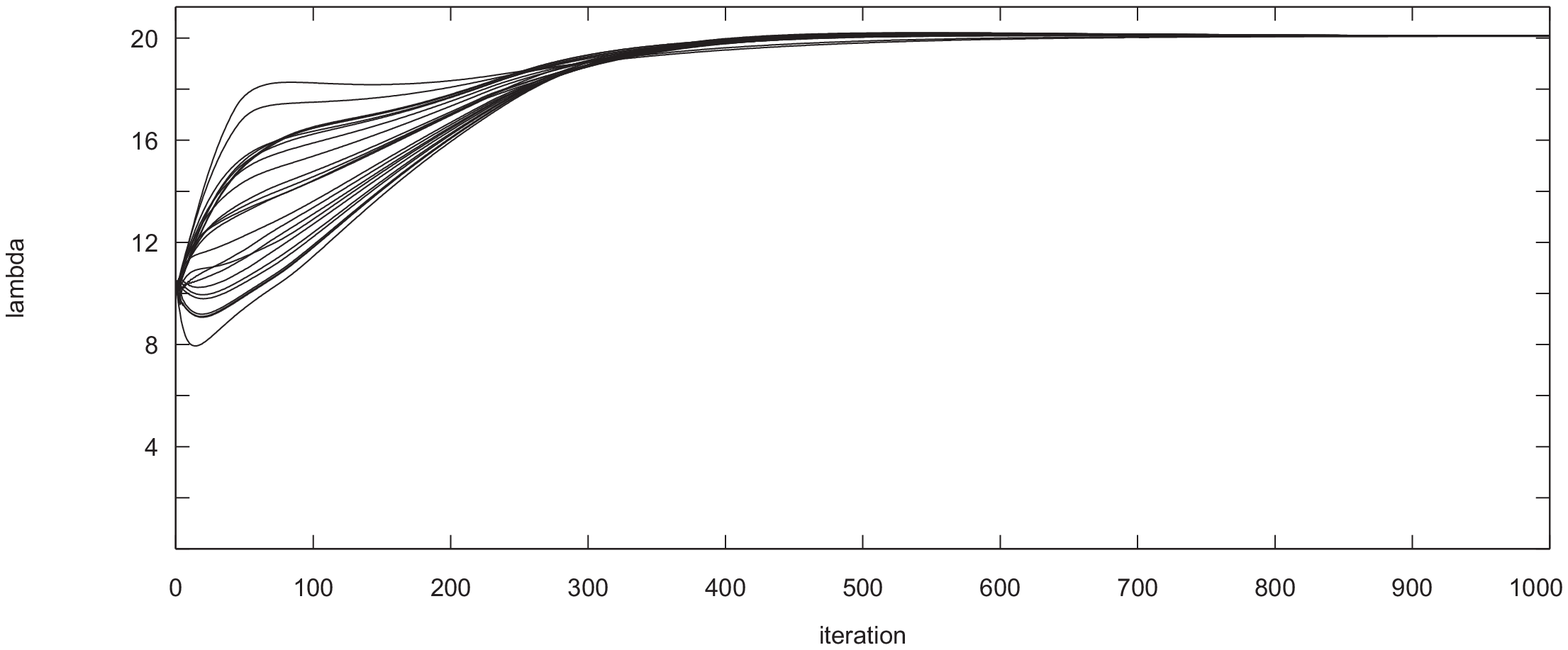}\label{lambdaSynch24}}
\caption{(a) Generation outputs, (b) bus angles, and (c) lagrange multipliers $\lambda$ (locational marginal prices) for uncongested case.}\label{SynchResults}
\end{figure}

In this first simulation, we keep the original line limits. Figure \ref{generationSynch24} gives the evolution of the total generation outputs for the buses to which generation is connected over the iterations, Fig.~\ref{thetaSynch24} the evolution of bus angles and Fig.~\ref{lambdaSynch24} the evolution of the Lagrange multipliers $\lambda$'s. One thousand iterations are displayed showing that convergence is achieved after about 600 iterations. The locational marginal prices all converge to the same value which is expected in an uncongested physical network and usage of DC power flow approximations.

\begin{figure}[t!]
\centering
 \footnotesize
 \psfrag{iteration}[c][c]{Iterations}
 \psfrag{obj}[c][c]{$|f-f^*|/f^*$}
 \psfrag{res}[c][c]{$\sum_i \sqrt{g_i^2}$}
 \psfrag{0}[c][c]{}
 \psfrag{100}[c][c]{}
 \psfrag{200}[c][c]{200}
 \psfrag{300}[c][c]{}
 \psfrag{400}[c][c]{400}
 \psfrag{500}[c][c]{}
 \psfrag{600}[c][c]{600}
 \psfrag{700}[c][c]{}
 \psfrag{800}[c][c]{800}
 \psfrag{900}[c][c]{}
 \psfrag{1000}[c][c]{1000}
 \psfrag{1}[r][r]{$10^1$}
 \psfrag{0}[r][r]{$10^0$}
 \psfrag{-1}[r][r]{$10^{-1}$}
 \psfrag{-2}[r][r]{$10^{-2}$}
 \psfrag{-3}[r][r]{$10^{-3}$}
 \psfrag{-4}[r][r]{$10^{-4}$}
    \subfigure[]{\includegraphics[width=8cm]{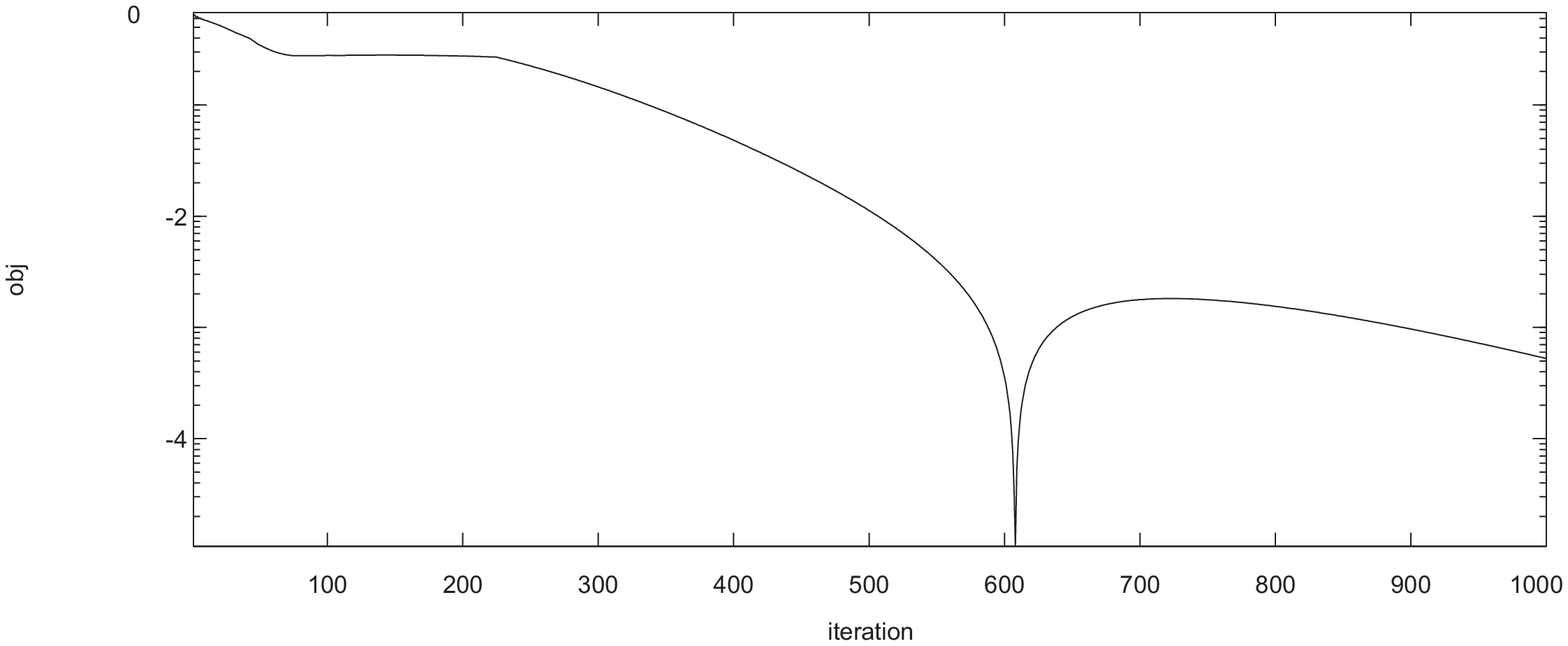}\label{logobjSynch24}}
    \subfigure[]{\includegraphics[width=8cm]{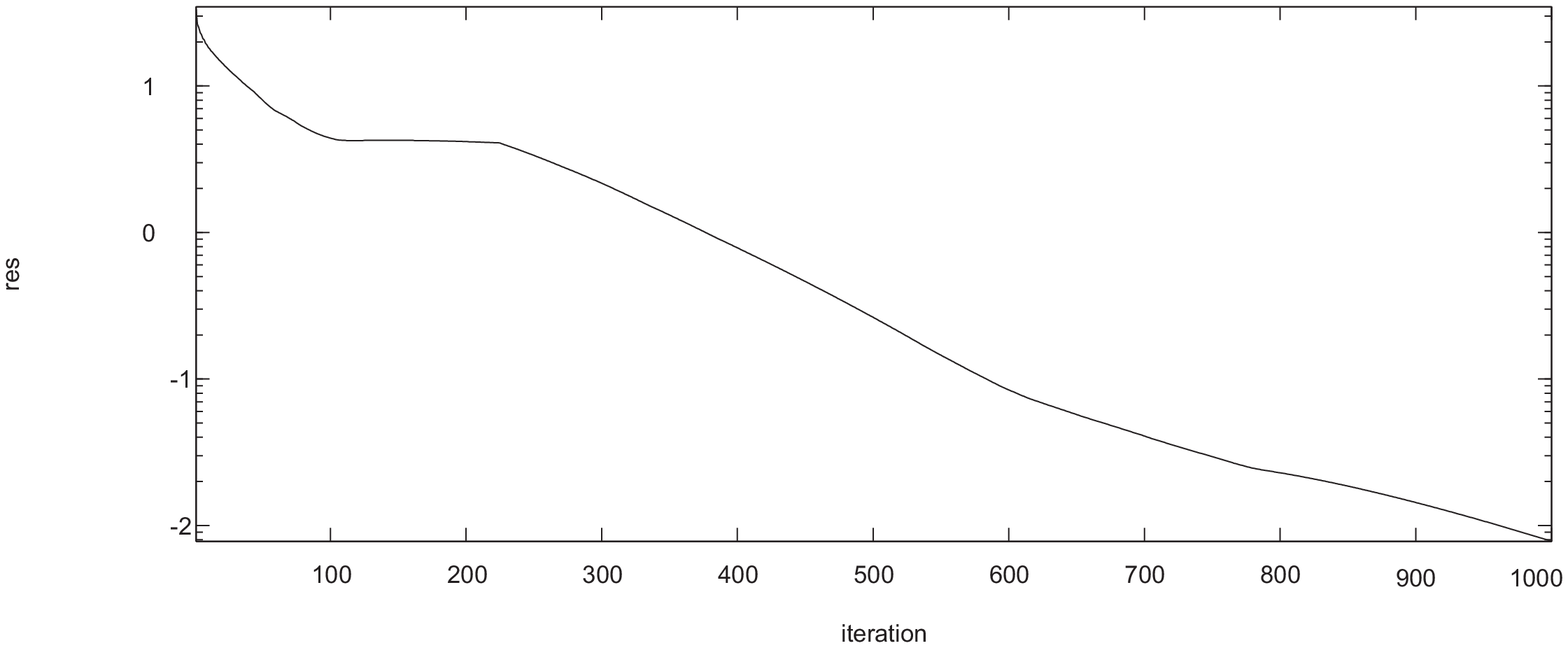}\label{logresSynch24}}
\caption{(a) Objective function value $\frac{|f-f^*|}{f^*}$ and (b) residual of equality constraints $\sum_i \sqrt{g_i^2}$ for uncongested case.}\label{logsynch}
\end{figure}

Figure \ref{logobjSynch24} provides the information on the relative distance, while \ref{logresSynch24} gives $rel$ the sum over the residuals of all power flow equations in the system. It is obvious that the error in the solution decreases fairly quickly as more iterations are carried out.

\vspace{0.1cm}
\begin{figure}[t!]
\centering
 \footnotesize
 \psfrag{iteration}[c][c]{Iterations}
 \psfrag{obj}[c][c]{$|f-f^*|/f^*$}
 \psfrag{res}[c][c]{$\sum_i \sqrt{g_i^2}$}
 \psfrag{0}[c][c]{}
 \psfrag{100}[c][c]{}
 \psfrag{200}[c][c]{200}
 \psfrag{300}[c][c]{}
 \psfrag{400}[c][c]{400}
 \psfrag{500}[c][c]{}
 \psfrag{600}[c][c]{600}
 \psfrag{700}[c][c]{}
 \psfrag{800}[c][c]{800}
 \psfrag{900}[c][c]{}
 \psfrag{1000}[c][c]{1000}
 \psfrag{1200}[c][c]{1200}
 \psfrag{1400}[c][c]{1400}
 \psfrag{1}[r][r]{$10^1$}
 \psfrag{0}[r][r]{$10^0$}
 \psfrag{-1}[r][r]{$10^{-1}$}
 \psfrag{-2}[r][r]{$10^{-2}$}
 \psfrag{-3}[r][r]{$10^{-3}$}
 \psfrag{-4}[r][r]{$10^{-4}$}
 \psfrag{-6}[r][r]{$10^{-6}$}
    \subfigure[]{\includegraphics[width=8cm]{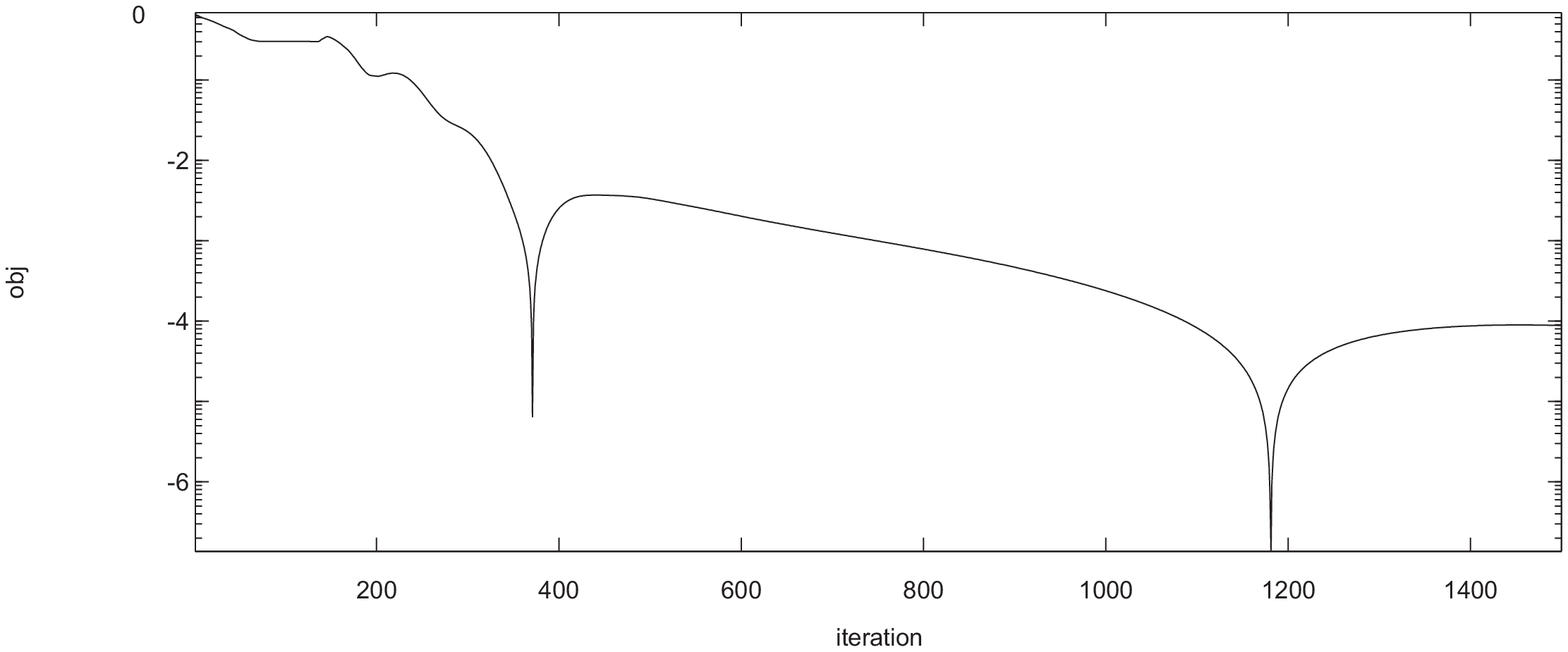}\label{logobjSynch24Con}}
    \subfigure[]{\includegraphics[width=8cm]{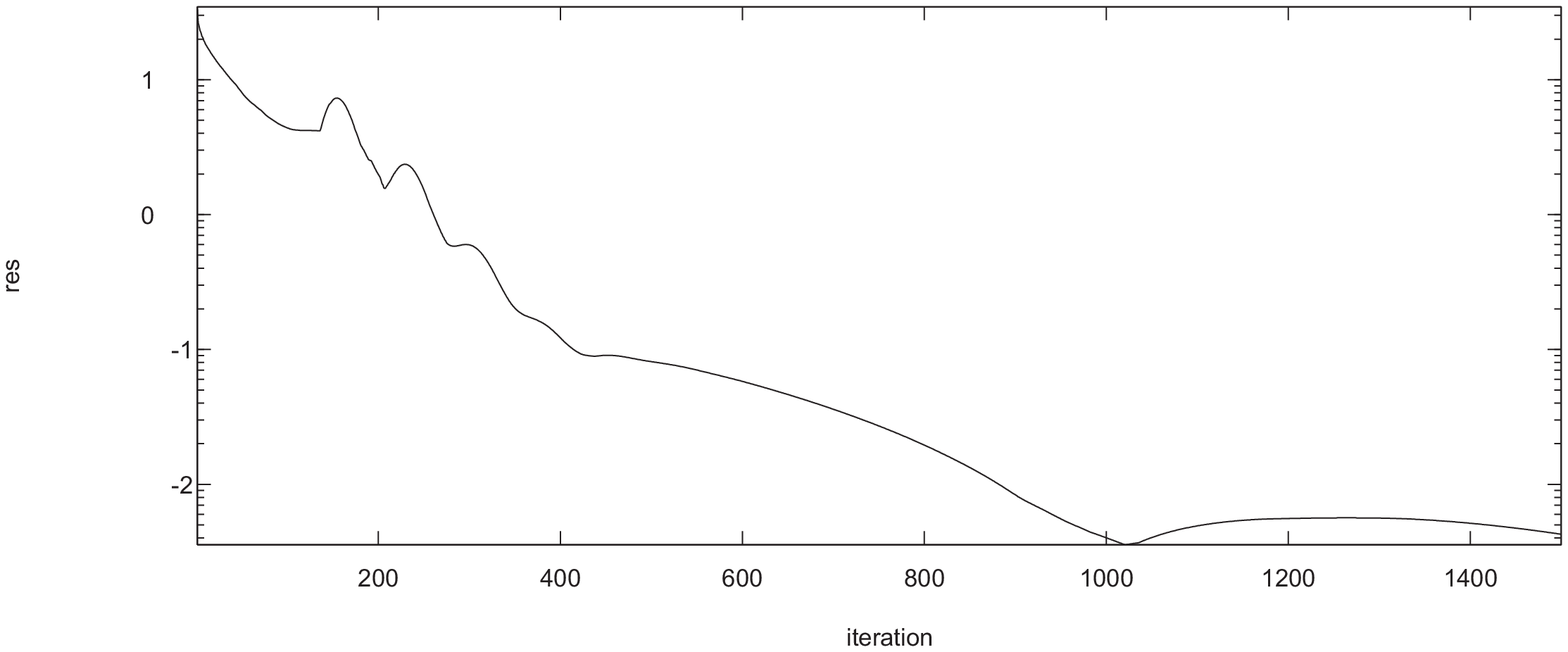}\label{logresSynch24Con}}
\caption{(a) Relative objective function value $\frac{|f-f^*|}{f^*}$ and (b) residual of equality constraints $\sum_i \sqrt{g_i^2}$ for congested case.}\label{logsynchCon}
\end{figure}

\subsubsection{Congested System}
\begin{figure}[t!]
\centering
 \footnotesize
 \psfrag{iteration}[c][c]{Iterations}
 \psfrag{theta}[c][c]{$\theta$ [rad]}
 \psfrag{generation}[c][c]{$P_{G}$ [MW]}
 \psfrag{lambda}[c][c]{$\lambda$}
 \psfrag{mulines}[c][c]{$\mu_{ij},~\mu_{ji}$}
 \psfrag{100}[c][c]{}
 \psfrag{200}[c][c]{200}
 \psfrag{300}[c][c]{}
 \psfrag{400}[c][c]{400}
 \psfrag{500}[c][c]{}
 \psfrag{600}[c][c]{600}
 \psfrag{700}[c][c]{}
 \psfrag{800}[c][c]{800}
 \psfrag{900}[c][c]{}
 \psfrag{1000}[c][c]{1000}
 \psfrag{1200}[c][c]{1200}
 \psfrag{1400}[c][c]{1400}
 \psfrag{5}[r][r]{5}
 \psfrag{10}[r][r]{10}
 \psfrag{15}[r][r]{15}
 \psfrag{20}[r][r]{20}
 \psfrag{25}[r][r]{25}
 \psfrag{0.2}[c][c]{0.2}
 \psfrag{0.1}[c][c]{}
 \psfrag{0.0}[c][c]{0.0}
 \psfrag{-0.1}[c][c]{}
 \psfrag{-0.2}[c][c]{-0.2}
 \psfrag{-0.3}[c][c]{}
 \psfrag{-0.4}[c][c]{-0.4}
 \psfrag{-0.5}[c][c]{}
 \psfrag{-0.6}[c][c]{-0.6}
 \psfrag{-0.7}[c][c]{}
 \psfrag{1}[r][r]{}
 \psfrag{2}[r][r]{200}
 \psfrag{3}[r][r]{}
 \psfrag{4}[r][r]{400}
 \psfrag{5}[r][r]{5}
 \psfrag{6}[r][r]{600}
    \subfigure[]{\includegraphics[width=8cm]{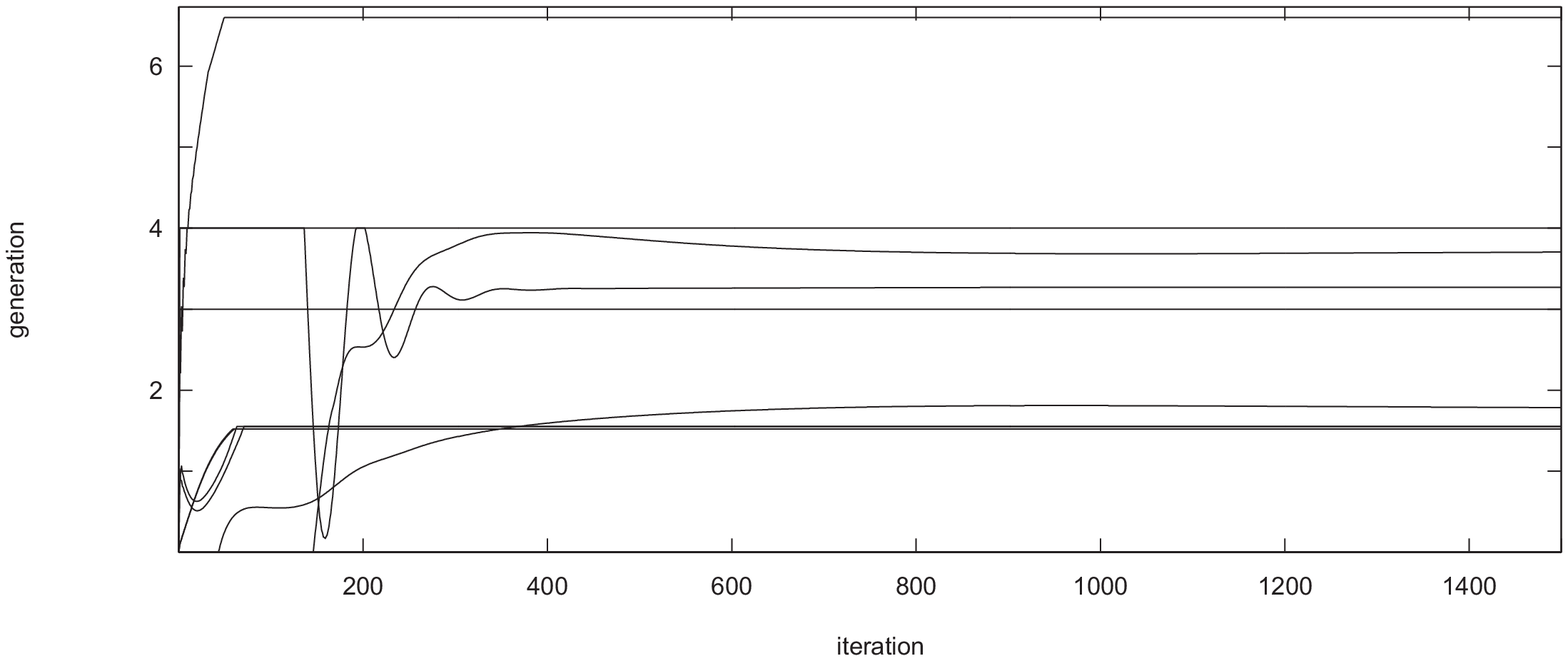}\label{generationSynch24Con}}
    \subfigure[]{\includegraphics[width=8cm]{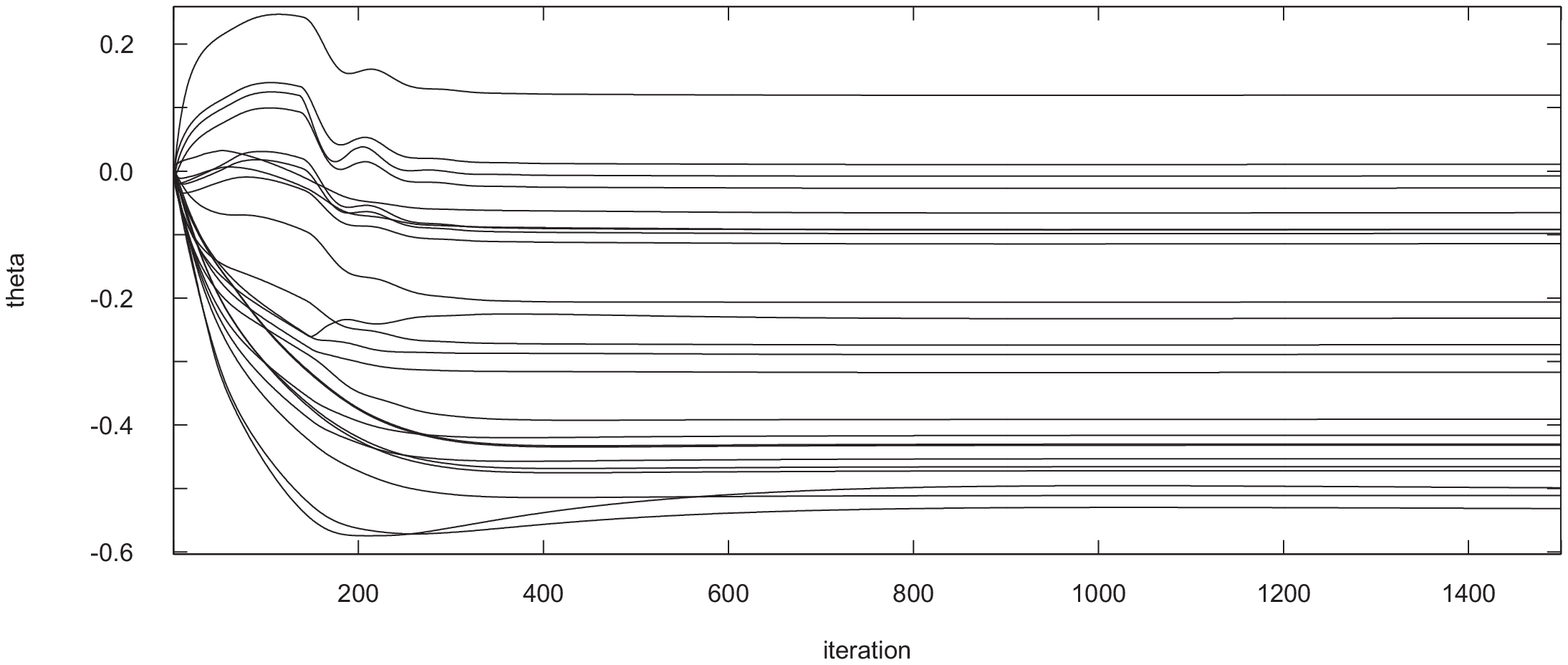}\label{thetaSynch24Con}}
    \psfrag{4}[r][r]{4}
    \subfigure[]{\includegraphics[width=8cm]{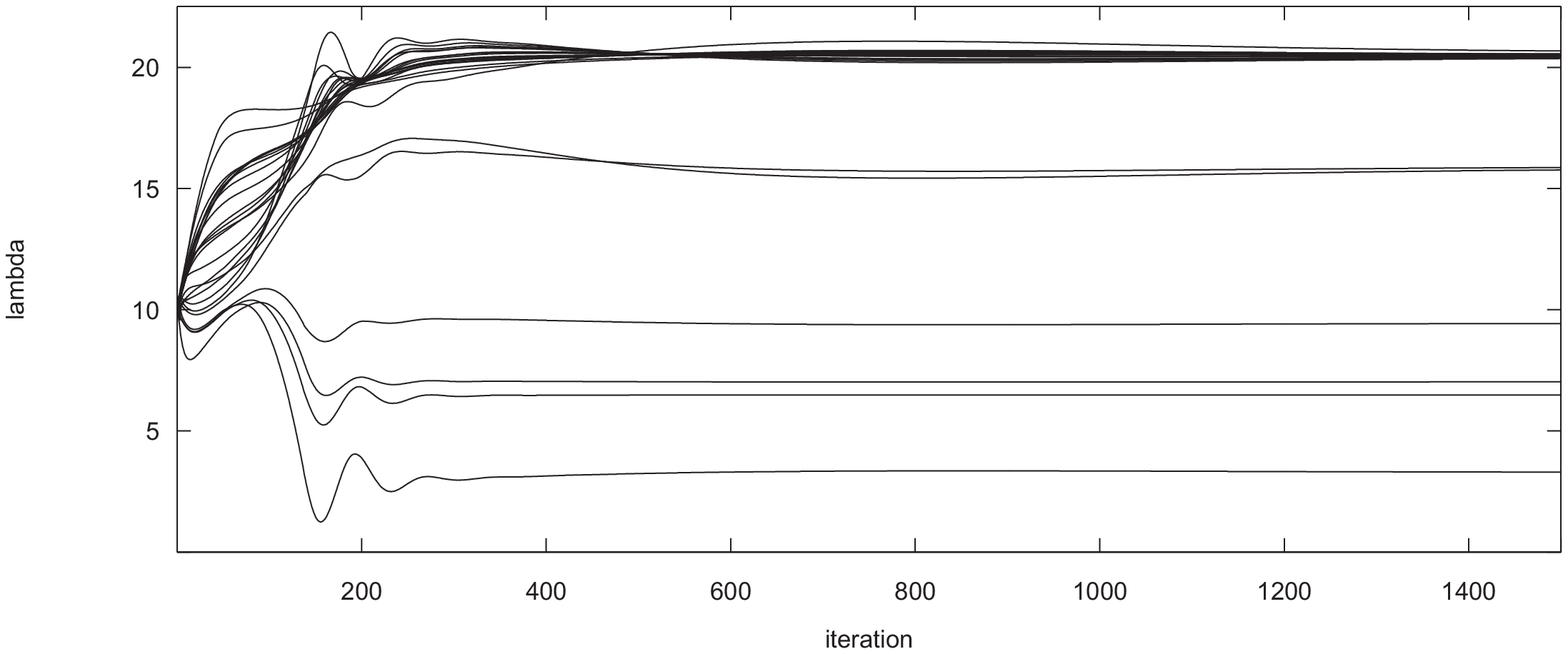}\label{lambdaSynch24Con}}
    \subfigure[]{\includegraphics[width=8cm]{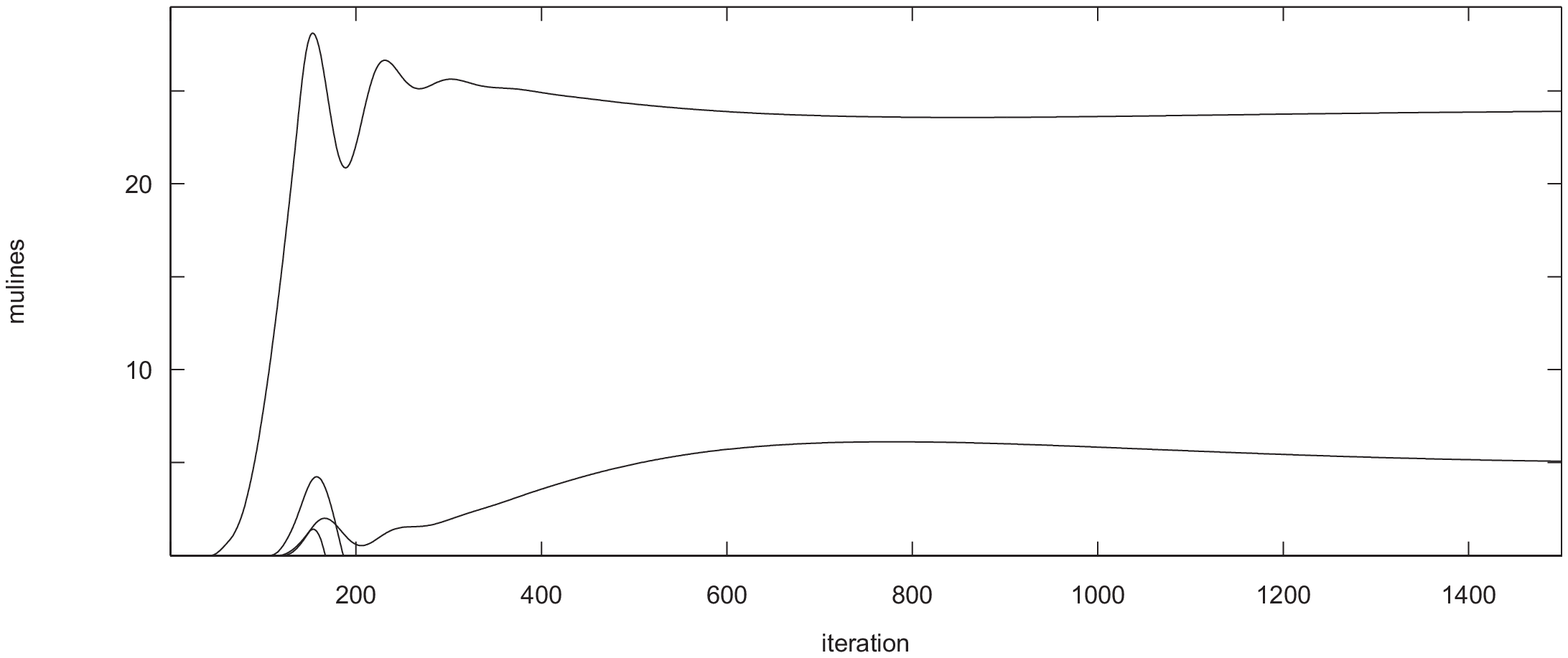}\label{mulinesSynch24Con}}
\caption{(a) Generation outputs, (b) bus angles, (c) lagrange multipliers $\lambda$ and (d) lagrange multipliers $\mu_{ij},~\mu_{ji}$ for congested case.}\label{SynchResultsCon}
\end{figure}

In order to create a situation in which lines reach their limits, all line limits are reduced to $55\%$ of their original values. This leads to the lagrange multipliers $\mu_{ij}$ associated with the line constraints for the congested lines to be non-zero for the positive flow direction. Hence, Fig.~\ref{SynchResultsCon} gives the values for the generation settings $P_G$, the bus angles $\theta$, the lagrange multipliers $\lambda$ and the Lagrange multipliers $\mu$ over the iterations for this case. Figure \ref{logsynchCon} provides the relative distance from the optimal cost function value and the sum of the constraint residuals.

Given that two lines reach their limits and the resulting non-zero update of the corresponding Lagrange multipliers, it takes more iterations to convergence than in the non-congested case. In addition, oscillations appear which could be prevented by reducing some of the tuning parameters but this would also lead to a larger number of iterations until convergence. It can be seen that the two lagrange multipliers associated with the line constraints of the congested lines are non-zero and the locational marginal prices $\lambda$ are not equal to the same value any more.

\section{Conclusion}\label{Conclusion}
In this paper, we presented a distributed approach to solve the DC Optimal Power Flow problem, i.e.~the generation dispatch is determined which minimizes the cost to supply the load in a distributed manner taking into account limited line capacities.
The main features of the algorithm are that it allows for a fully distributed implementation down to the bus level without the need for a coordinating central entity, the individual updates per iteration consist of simple function evaluations and exchange of information is limited to bus angles and Lagrange multipliers associated with power flow equations and line constraints among the neighboring buses.
In particular, there is no need to share information about generation cost parameters or generation settings.

The algorithm was tested in the IEEE Reliability Test System showing that it converges to the overall optimal solution. Moreover, this paper discusses the convergence criteria for the proposed distributed method, and analytically proves that the limit point of our innovation-based approach is the optimal solution of the OPF problem.

\section*{Acknowledgment}
The authors would like to thank ARPA-E and SmartWireGrid for the financial support for this work under the project ``Distributed Power Flow Control using Smart Wires for Energy Routing'' in the GENI program.

\appendix
\subsection*{A Sufficient Condition for the Convergence of the Proposed Algorithm}\label{Conv}
The following section provides a sufficient condition for the convergence of the proposed distributed algorithm given by (\ref{iterative0}). To this end, the following assumption on the matrix $A$ as defined in (\ref{updateMatrix}) is imposed:\\
\noindent \textit{A.1:} There exists an $\ell_{p}$-norm such that the tuning parameters $\alpha$, $\beta$, $\gamma$ and $\delta$ can be designed to achieve $\left \| I-A \right \|_{p}<1$.\\
\indent \textit{Remark~1:} Note that the projection operator $\mathbb{P}$ in our context involves component-wise projections, and, hence, non-expansive with respect to $\ell_{p}$-norms, i.e., the following holds for any two iterates $X(k)$ and $X(v)$,
 \begin{equation}
 \left \|\mathbb{P}(X(k))-\mathbb{P}(X(\nu))  \right \|_{p}\leq \left \|X(k)-X(\nu)  \right \|_{p}.\label{contraction}
 \end{equation}
 \indent \textit{Remark~2:} Based on \textit{Remark~1}, the following equations hold,
\begin{align}\label{remark3}
& \left \| \widetilde{X}(k+1)- \widetilde{X}(k)\right \|_{p} \\ \nonumber
&=\left \| \mathbb{P}[(I-A)\widetilde{X}(k)+C]-\mathbb{P}[(I-A)\widetilde{X}(k-1)+C] \right \|_{p}\\ \nonumber
&\leq\left \| (I-A)\widetilde{X}(k)+C-(I-A)\widetilde{X}(k-1)-C \right \|_{p}\\ \nonumber
&\leq\left \| (I-A)\right \|_{p} \left \| \widetilde{X}(k)-\widetilde{X}(k-1)\right \|_{p}
\end{align}
Consequently, (\ref{remark3}) leads to
\begin{equation}\label{DiffX}
\left \| \widetilde{X}(k+1)-\widetilde{X}(k)\right \|_{p}\leq{\left \| I-A\right \|}^{k}_{p} \left \| \widetilde{X}(1)-\widetilde{X}(0)\right \|_{p}
\end{equation}
\indent \textit{Theorem~3:} Let \textit{A.1} hold, then the algorithm presented in (\ref{iterative0}) achieves convergence.
 %i.e. $\lim_{k\rightarrow \infty} \widetilde{X}(k)=\widehat{X}\Rightarrow \\
% ~~~~~~~~\forall~ \delta \geq 0~~ \exists ~K_{\delta}~~s.t.~\left | \widetilde{X}(k)-\widehat{X} \right |<\delta~~\forall~k \geq K_{\delta}$
\\ \indent \textit{Proof:} The distance between the values of $\widetilde{X}$ at two iterations $k$ and $\nu$ is given by,
\begin{align}\label{DiffX1}
&\left \| \widetilde{X}(k)-\widetilde{X}(\nu)\right \|_{p}=\\ \nonumber
&\left \| \widetilde{X}(k)-\widetilde{X}(k-1)+\widetilde{X}(k-1)- \cdots +\widetilde{X}(\nu+1)-\widetilde{X}(\nu)\right \|_{p} \\ \nonumber
&\leq \left \| \widetilde{X}(k)-\widetilde{X}(k-1) \right \|_{p}+ \cdots +\left \|  \widetilde{X}(\nu+1)-\widetilde{X}(\nu)\right \|_{p}.\\ \nonumber
\end{align}
Moreover, using (\ref{DiffX}) the following equation can be derived:
\begin{align}\label{DiffX2}
&\left \| \widetilde{X}(k)-\widetilde{X}(\nu)\right \|_{p}\\ \nonumber
&\leq \left \| \widetilde{X}(k)-\widetilde{X}(k-1) \right \|_{p}+ \cdots +\left \|  \widetilde{X}(\nu+1)-\widetilde{X}(\nu)\right \|_{p}\\ \nonumber
&\leq\left (  {\left \| I-A\right \|}^{k}_{p}+\cdots +{\left \| I-A\right \|}^{\nu}_{p}\right ) \left \| \widetilde{X}(1)-\widetilde{X}(0)\right \|_{p}.
\end{align}
Since \textit{A.1} holds, we have
\begin{equation}
\lim_{\nu\rightarrow \infty}\|I-A\|^{\nu}_{p}=0\label{I-Ainf}.
\end{equation}
Hence, combining (\ref{DiffX2}) and (\ref{I-Ainf}) further implies,
%\begin{equation}
%\forall~ \epsilon>0,~\exists~N~s.t~~ k,\nu>N~\Rightarrow \left \| X^{k}-X^{\nu} \right \|\leq\epsilon\label{CasuchySeq}
%\end{equation}
%And because of the property of projection, (\ref{contraction1}), the following equation holds,
\begin{equation}
\forall~ \epsilon>0,~\exists~N~s.t~~ k,\nu>N~\Rightarrow \left \| \widetilde{X}(k)-\widetilde{X}(\nu) \right \|_{p}\leq\epsilon\label{CasuchySeq}
\end{equation}
Therefore the sequence of $\left \{ \widetilde{X}(i) \right \}_{i=0}^{\infty}$, which is introduced by (\ref{iterative0}), is a Cauchy sequence. Since a sequence of real vectors converges to a limit in $\mathbb{R}^n$ if and only if it is Cauchy, it follows that the proposed iterative algorithm is convergent, i.e., $X(i)\rightarrow X^{\ast}$ as $i\rightarrow\infty$ for some $X^{\ast}\in\mathbb{R}^{n}$.
\vspace{0.05cm}
\bibliographystyle{IEEE}
\bibliography{References1}

\vfill

\end{document}